\documentclass[12pt]{article}
\usepackage{amssymb,amsbsy,amsmath,amsfonts,amssymb,amscd}
\usepackage{latexsym}

\font\Bbb=msbm12

\begin{document}
\def\DD{\Delta}
\def\dd{\delta}
\def\oo{\omega}
\def\ee{\epsilon}
\def\OO{\Omega}
\def\phi{\varphi}
\font\SY=msam10
\def\square{\hbox{\SY\char03}}
\def\wt{\widetilde}
\def\a{\alpha}
\def\b{\beta}
\def\Deg{{\rm Deg}\,}
\def\deg{{\rm deg}}
\def\ggg{\delta}
\def\res{{\rm res}}
\def\p{\rho}
\def\gg{\gamma}
\font\goth=eufm10 scaled1200
\def\sm{\hbox{\goth S}}
\def\Gl{{\rm GL}\,}
\def\M{{\rm M}}
\def\BZ{{\mathbb Z}}
\def\BP{{\mathbb P}}
\def\BF{{\mathbb F}}
\def\BG{{{\mathbb G}}}
\def\BN{{\mathbb N}}
\def\BR{{{\mathbb R}}}
\def\BC{{{\mathbb C}}}
\def\BQ{{{\mathbb Q}}}
\def\BT{{\mathbb T}}

\newtheorem{lemma}{Lemma}[section]
\newtheorem{theorem}[lemma]{Theorem}
\newtheorem{corollary}[lemma]{Corollary}
\newtheorem{proposition}[lemma]{Proposition}

\title{The conjugacy problem for two-by-two matrices over polynomial rings}
\author
{Fritz J. Grunewald$^{1}$,  Natalia K. Iyudu$^{2}$}

\date{}

\maketitle

{\small\tt\noindent \llap{$^1\,$}Mathematisches Institut, Henrich
Heine Universit\"at, \hfill\break\noindent 40225 D\"usseldorf,
Germany,
\\
\llap{$^2\,$}Department of Pure Mathematics, Queen's University Belfast,
 \hfill\break\noindent Belfast BT7 1NN, U.K.}

\medskip

{\bf e-mails:\ \tt n.iyudu@qmul.ac.uk,\ \
fritz@math.uni-duesseldorf.de\\}

\maketitle
\bigskip
\bigskip
\bigskip

\begin{abstract}
We give an effective solution of the conjugacy problem for
two by two matrices over the polynomial ring in one variable over a finite
field.
\end{abstract}

\hrule
\bigskip
\bigskip
\bigskip
\tableofcontents
\medskip

\bigskip
\bigskip

\section{Introduction}

We consider here the conjugacy problem in the ring of two by two
matrices $M(2,\mathbb F [x])$ over the polynomial ring $\BF[x]$,
where $\BF$ is a finite field. We say that two matrices $A,\, B\in
\M_2(\BF[x])$ are conjugate if there is a conjugating matrix $U$ in
the group $\Gl(2,\BF[x])$ of invertible matrices over $\BF[x]$, such
that $U$ satisfies $B=UAU^{-1}$. In the following we write $\deg(p)$
for the degree of a polynomial $p \in \BF[x]$ and $\deg(A)$ for the
maximal degree of the entries of $A\in \M_2(\BF[x])$. We prove:

\begin{theorem}\label{teoa}
Let $\BF$ be a finite field with $q$ elements and $A,\, B\in
\M_2(\BF[x])$. Let $\dd$ be the maximum of $\deg(A),\, \deg(B)$. If
$A,B$ are conjugate, then there is a conjugating matrix $U$ with
$\deg(U)\leq (1+q) \dd q^{7 \dd} $.
\end{theorem}

For certain pairs of matrices $A,\, B\in \M_2(\BF[x])$ the estimate
of the degrees of the entries of $U$ can be improved to be linear in
$\delta$ not depending on $q$ (see  Proposition \ref{Proposition
4.1.}). Theorem \ref{teoa} shows that there is an algorithm which
decides whether two matrices $A,\, B\in \M_2(\BF[x])$ are conjugate
or not. Hence we can state:

\begin{corollary}\label{coroa}
Let $\BF$ be a finite field, then the conjugacy problem in the group
$\Gl(2,\BF[x])$ is effectively solvable.
\end{corollary}

Corollary \ref{coroa} should be compared with the solution of the
conjugacy problem in an arithmetic group. The conjugacy problem for
$\Gl(n,\BZ)$ ($n\in\BN$) was solved in \cite{FG}. But even in the
case $n=2$ no explicit estimates like those from Theorem \ref{teoa}
are known. Also the algorithms described in \cite{GS}, which solve
the conjugacy problem in any arithmetic group, do not give estimates
for the degree of a conjugating matrix. The method of solution
employed in \cite{FG} for the case of $\Gl(n,\BZ)$ ($n\in\BN$) can
be extended (without giving any estimates) to the case of
$\Gl(n,\BF[x])$ ($n\in\BN$, $\BF$ a finite field) when the
characteristic of the field $\BF$ does not divide the size $n$ of
the matrices. Also our method in \cite{FG} provides extra
difficulties in case $\BF$ has characteristic 2. Further features of
the conjugacy problem in $\Gl(2,\BF[x])$ are described in Section
\ref{finsec}.

Given a matrix $A\in \Gl(2,\BF[x])$ we define
\begin{equation}\label{centra}
{\rm Z}(A):=\{\, U \in \Gl(2,\BF[x])\ \mid\  U A U ^{-1}=A\, \}
\end{equation}
to be its centralizer. In case $A\ne {\bf 1}$ is semisimple it is
well known that ${\rm Z}(A)$ is either finite or the direct product
of an infinite cyclic group by a finite group. By our methods we can
give an estimate for the degrees of the entries of a generator of
the infinite part:

\begin{theorem}\label{teob}
Let $\BF$ be a finite field with $q$ elements and $A\in \Gl(2,
\BF[x])$ a semisimple matrix, not equal to the identity matrix, such
that Z($A$) is infinite. Then there is a matrix $U\in {\rm Z}(A)$
which generates Z($A$) up to a finite group with $\deg(U)\le \deg(A)
q^{2\deg(A)}$.
\end{theorem}

Our method to prove Theorem \ref{teoa} uses a reduction to a
quadratic equation in two variables. As a special case Pell's
equation
\begin{equation}\label{pell}
u^2+Dv^2=1
\end{equation}
with $D\in\BF[x]$ arises. Let us call $D\in\BF[x]$ to be positive if
it is neither constant nor a square, has even degree and highest
coefficient a square. Let us furthermore call a solution $(u,v)$ of
(\ref{pell}) trivial if $u,\, v\in \BF$ holds. We prove:

\begin{theorem}\label{teoc}
Let $\BF$ be a finite field with $q$ elements and $D\in \BF[x]$ a
positive polynomial. Then (\ref{pell}) has a nontrivial solution
$(u,v)$ with $\deg(u),\, \deg(v)\le q^{\deg(D)}$.
\end{theorem}

Pell's equation (\ref{pell}) has been studied extensively in the
paper of Emil Artin of 1924 \cite{EA}. He investigate Pell's
equation through continued fraction expansions. But he assumed that
the characteristic of $\BF$ is not equal to 2. Our result in Theorem
\ref{teoc} follows straightforward from [1]. We then modify Artin's
technique for the case of characteristic $2$.

Following our reduction we have to analyze the solution of the
general quadratic equation
$$au^2+buv+cv^2=d$$
where $a,\dots,d$ are polynomials in $\BF[x]$. We use new degree
function on certain quadratic extension rings of $\BF[x]$ in the
imaginary case and control the behavior of continued fraction
expansion in the real case.

Note that as we will see while the case of characteristic 2 in some
sense more difficult the estimations for algorithms in this case
which we obtain turned out be better than in the case of the other
finite fields.

It worth mentioning also that the reduction step itself gives quite
noticeable impact to the whole estimate of the degree of conjugating
matrix.

\section{Reduction to a quadratic equation}\label{redu}

 Let $\BF$ be a field. We consider here pairs of matrices
$$A=\left( \begin {array}{cc}
                a_{11}&a_{12}\\
                a_{21}&a_{22}
        \end {array} \right),
B=\left( \begin {array}{cc}
                b_{11}&b_{12}\\
                b_{21}&b_{22}
        \end {array} \right)\in M_2(\BF[x]),$$
which we call rationally conjugate if they are conjugate by an
element of $\Gl(2,\BF(x))$, where $\BF(x)$ is a field of rational
functions over $\BF$. Being rationally conjugate implies the
conditions:
$$Tr(A)=Tr(B), det(A)=det(B).$$
Suppose we want to find a matrix
$$U=\left( \begin {array}{cc}
                u&p\\
                v&q
        \end {array} \right) \in GL_2(\BF[x]),$$
which conjugates $A$ to $B$. We are lead then to four linear
equations given by the matrix entries of $UA-BU$ in the variables
$u,p,v,q,$ plus the quadratic equation $det(U) \in \BF^*=\BF
\backslash \{0\}$. Elementary considerations of this system of
equations proves:

\begin{lemma}\label{ut}
Suppose the matrices $A,B \in \M_2(\BF[x])$ satisfy
$a_{21}=b_{21}=0$. Let $\delta$ be the maximum of the degrees of the
entries of $A,B$. Then $A,B$ are conjugate by an element of $\Gl(2,
\BF[x])$ if and only if they are conjugate by $U\in \Gl(2,\BF[x])$
with $deg(U) \leq \delta$.
\end{lemma}

The above lemma proves Theorem \ref{teoa} in the special case of
upper triangular matrices $A,B$.

Let $char \, \BF = 2$. The conjugating condition $UAU^{-1} = B$ for
$U \in \Gl(2,\BF[x])$ is equivalent to the system
\begin{equation}
                uq+pv \in \BF^*
\label{uqpv*}
\end{equation}
\begin{equation}
                UA=BU
\label{uabu}
\end{equation}
where $ U=\left( \begin {array} {cc}
                u&p\\
                v&q
        \end {array} \right)$,
$\BF^*=\BF \backslash \{0\}$ is the multiplicative group of the
field $\BF$. Since the set of conjugating matrices is stable under
multiplication by a non-zero constant and we have a unique square
root in our field, the solvability of the system (\ref{uqpv*},
\ref{uabu}) is equivalent to the solvability of the same system with
(\ref{uqpv*}) replaced  by
\begin{equation}
                uq+pv=1
\label{uqpv}
\end{equation}

The quadratic equation (\ref{uqpv})
with additional linear conditions  (\ref{uabu})
can be reduced
to one quadratic equation in two different ways.

First, using the procedure of construction of the generating system
of syzygies module (\cite{AdLo}) for the linear system (\ref{uabu}).
Another way is a direct substitution of the solution in rational
functions of the linear system (\ref{uabu}). In both cases we obtain
an equation of the type $au^2+buv+cv^2=d$, but in the first case
with variables of different meaning,  in the second case with some
additional divisibility conditions. We will follow the second way.
Let $A=\left( \begin {array}{cc}
                a_{11}&a_{12}\\
                a_{21}&a_{22}
        \end {array} \right)$,
$B=\left( \begin {array}{cc}
                b_{11}&b_{12}\\
                b_{21}&b_{22}
        \end {array} \right)$  $\in M_2(\BF[x]) $,
substituting rational solutions
\begin{equation}
p= \frac{(a_{11}+b_{11})u+b_{12} v}{a_{21}},
 q=\frac{b_{21} u + (a_{11}+b_{11}) v}{a_{21}}
\label{pq}
\end{equation}
from (\ref{uabu}) to  (\ref{uqpv}) we obtain the equation

\begin{equation}
  b_{21} u^2+(b_{11} + b_{22})uv +  b_{12} v^2 =  a_{21}.
\label{Q}
\end{equation}

To ensure that $p$ and $q$ are polynomials, for the solutions $u,v$
of (\ref{Q}) we have to check that $a_{21} \Big|
(a_{11}+b_{11})u+b_{12} v$ and    $a_{21} \Big| b_{21} u + (a_{11}+
b_{11}) v$.

\vspace{5mm}

We can consider separately quite an easy  case when one of matrices
is diagonal. Matrices $A,B$ where B is diagonal are conjugate if and
only if $\frac{a_{12}(b_{22}-b_{11})} { g.c.d.(a_{11}-b_{11},
a_{12}) g.c.d.(a_{11}-b_{22}, a_{12})} \in \BF^*$ (in case when
$(a_{11}-b_{11}, a_{12}) \neq (0,0)$ and $(a_{11}-b_{22}, a_{12})
\neq (0,0)$). In case $(a_{11}-b_{11}, a_{12})= (0,0)$ the condition
looks slightly different:
$\frac{(b_{11}-a_{22})(b_{22}-a_{11})-a_{12}a_{21}}
{g.c.d.(b_{11}-a_{22}, a_{12}) g.c.d.(a_{11}-b_{22}, a_{12})} \in
\BF^*.$

  After that we can restrict ourselves by the case
$(a_{12}, a_{21}) \neq (0,0)$ and
$(b_{12}, b_{21}) \neq (0,0)$.
Since any matrix is conjugate with it's transposed,
we can assume without loss of generality that
$a_{21}$ and $b_{21}$ are non-zero.

\section{Preliminary considerations for the solution of
the quadratic equation $au^2+buv+cv^2=d$}\label{prelim}

We study here the equation $au^2+buv+cv^2=d$ with $a \neq 0$.
Multiplying the equation  by $a$ and making the change of variables
$u_1=au, \, v_1=v$ we obtain the equation
$u_1^2+bau_1v_1+cav_1^2=da$ with the monic polynomial in the
lefthand side. After we find a solution, we have to check whether
$a$ divides $u_1$ and only in this case $a|u_1$ will give a solution
of the initial equation.

\vspace{5mm}

We consider the following three cases determined by the nature of
the roots of the  equation
\begin{equation}
t^2+bt+c=0.  \label{t}
\end{equation}

Case 1. Equation (\ref{t}) is solvable in $\BF(x)$. Note, that in
fact it means that (\ref{t}) is solvable in  $\BF[x]$. If
$\frac{F}{Q}$ is a rational solution and ${\rm g.c.d.}(F,Q)=1$, then
from $F^2+bFQ+cQ^2=0$ follows $Q \Big| F^2$, hence $Q$ can be only
constant and the solution is in fact polynomial. In this case
$d=u^2+buv+cv^2=(u+ \Delta v)(u+ (b+\Delta) v)$ is a product of two
polynomials. There exists a finite set of factorizations of $d$ into
two multiples from $\BF[x]$, hence we obtain a finite number of
linear systems on $u,v$.

Let us mention that if $b = 0$ then we have a rational solutions (we
are in case 1). Indeed, let
 $b=0$:
$u^2+cv^2=d$. We can just present the coefficients as follows:
$c(x)=c_0(x^2)+xc_1(x^2), d(x)=d_0(x^2)+xd_1(x^2)$. If
$u=\sum\limits_{i=0}^n u_i x^i$, let $\widetilde
u=\sum\limits_{i=0}^n u_i^2 x^i$, then $u^2(x)=\widetilde u (x^2)$.
Considering separately cases of even and odd degrees on $x$ we get
two linear equations on $\widetilde u $ and $\widetilde v $:
$$\left\{ \begin {array} {l}
   \widetilde u  + c_0 \widetilde v =d_0, \\
   c_1 \widetilde v =d_1.
        \end {array}\right.
$$
If $c_1=d_1=0$ there are infinitely many rational solutions: $\wt
u=-\frac{d_0}{c_0}\wt v$. Otherwise, there is at most one rational
solution. For any solution $\widetilde u$ and $\widetilde v $ we can
uniquely determine $u$ and $v$: $u_i = \sqrt{\widetilde u_i}$ and
$v_i = \sqrt{\widetilde v_i}$, --- due to the existence and
uniqueness of square roots in our field $\BF$.

The other two cases are more essential and will compose our main
treatment later on.

From this point we will suppose that the equation (\ref{t}) has no
rational solutions, particularly  $b \neq 0.$

Let us define a {\it degree function on} $\BF[x]$ as an ordinary
degree on non-zero polynomials and ${\rm deg\,(0)=-\infty}$.

We deal with the completion of  $\BF[x]$ by the valuation $\vert p
\vert =2^n$, where $n= {\rm deg\,}\, p$ (valuation of zero is 1).
This completion is the algebra of formal power series $
K=\BF((x))=\{\sum\limits_{i=-\infty}^d \alpha_i x^i, \alpha_i \in
\BF, d \in {\BZ}$\}.

Here is an essential in what follows

\noindent {\bf Definition. }
{\it Let $\p=\sum\limits_{-\infty}^d \a_{n}x^{n}$
be a power series. We say that $d$ is a degree of $\p$ if
$\a_d \neq 0$.}

Case 2. If (\ref{t}) is solvable in $K \backslash \BF[x]$, we  say
that it is a {\it real} case.

Case 3. In case when the solution can not be presented
as a power series,
we say it is an {\it imaginary} case.

\vspace{2mm}

Let us consider now the ring (and corresponding function field)
$R=\BF[x,t]/f_x(t)$, where $f_x(t)= t^2+bt+c$, $b,c \in \BF[x]$.
Obviously elements of $R$ can be uniquely presented as $u+\Delta v$,
where $u,v \in \BF[x]$ and $f_x(\Delta)=0$. We can define the {\it
norm} of an element $\omega = u+ \Delta v$ as $N(u+ \Delta
v)=F(u,v)=u^2+buv+cv^2$. Let define the {\it conjugate} element for
$\oo$ as follows: $\oo'= u+ (b+\Delta) v$.

It is easy to check that the introduced notions of norm and
conjugate element satisfy  the natural properties.


\begin{lemma}
\label{Lemma 2.1.}
a). $N(\oo)=\oo \oo'$; \
b). $(\oo_1 \oo_2)'= \oo_1' \oo_2'$.\
c). $N(\oo_1 \oo_2)=N(\oo_1)N(\oo_2).$ \
d). $N(\oo^{-1})=N(\oo)^{-1}.$
\end{lemma}


\begin{lemma}
\label{Lemma 2.5.} An element $\epsilon$ is a unit of $R$ if and
only if $N(\epsilon) \in \BF^*.$
\end{lemma}
{\bf Proof}. If ${\ee}^{-1}$ does exist then by lemma \ref{Lemma
2.1.} c). and d). $N(\ee \ee^{-1})=N(\ee)N(\ee^{-1})=
N(\ee)N(\ee)^{-1}= 1$.  Hence $N(\ee)$ is an invertible polynomial,
i.e. $N(\ee) \in \BF^*$.  \square

We treat real and  imaginary  cases in different ways, thus we need
first to be able to distinguish these cases. The following
proposition serve for this.

Our equation as earlier is
\begin{equation}
u^2+buv+cv^2=d.  \label{stard}
\end{equation}

\begin{proposition}\label{Theorem 3.1.}\ \break

{\rm I}. If $\rm deg\, c>2 \rm deg\, b$ and $ {\rm deg\,}  c$ is odd
then we
are in the imaginary case.\\

{\rm II}.  If $ {\rm deg\,}  c\ge 2 {\rm deg\,} b$ and $ {\rm deg\,}
c$ is even then there exists an invertible linear change of
variables which turns the equation into new one $u^2+buv+\widetilde
cv^2=d$,
with $ {\rm deg\,}  \widetilde c< {\rm deg\,}  c$.\\

{\rm III}.  If $ {\rm deg\,}  c = 2 {\rm deg\,} b$, we have two
possibilities. In case if the equation  $ b_0^2 t^2 + b_0^2 t
+c_0=0$ ( $b_0, c_0 \in \BF$ -- coefficients near highest terms of
$b $ and $c$ ) is solvable in $\BF$, there exists a change of
variables which turns the equation into new one $u^2+buv+\widetilde
cv^2=d$, with $ {\rm deg\,} \widetilde c< {\rm deg\,} c$.
Otherwise  we are in the imaginary case.\\

{\rm IV}. If $ {\rm deg\,}  c< 2 {\rm deg\,} b$ we are in the real case.\\
\end{proposition}
{\bf Proof.} I. Suppose that $\oo = \sum\limits_{- \infty}^m a_n t^n
\in K$ is a root of the equation $\oo^2+b\oo+c=0$, $a_n \neq 0$. It
is necessary for the cancellation that the degrees of a pair of
terms in this equation are equal and the degree of the third one is
not grater than that. A priori there exist three possibilities.

1). ${\rm deg\,} \oo^2 = {\rm deg\,} b \oo$. Then ${\rm deg\,}(
 \oo ^2+b\oo) \leq 2 {\rm deg\,}b < {\rm deg\,} c$ and cancellation is
in fact impossible.

2). ${\rm deg\,}b \oo = {\rm deg\,} c$. Hence $ {\rm deg\,} \oo =
{\rm deg\,}c - {\rm deg\,}b$ and ${\rm deg\,}(b \oo + c) \leq {\rm
deg\,}c < {\rm deg\,} \oo^2$ and cancellation is again impossible.

3). ${\rm deg\,} \oo^2 = {\rm deg\,} c$. This case is not possible because
 ${\rm deg\,} c$ have to be odd.

This means that there are no solutions of (\ref{stard}) in power
series.

II.

We try to find desired change of variables in the form:
$$\left\{\begin{array}{l}
u'=u+brv,\\
v'=v\end{array}\right.
$$
In that case we have to find $r$ such that ${\rm deg\,}
(b^2r^2+b^2r+c) < {\rm deg\,} c$. Highest terms in this sum are
$b^2r^2$ and $c$ ($\deg c > 2 \deg b,$ hence $\deg r \neq 0$). To
provide their cancellation we take ${\rm deg\,} r = \frac{1}{2}
({\rm deg\,}c-2{\rm deg\,}b)$. Coefficients near the highest terms
also have to coincide: $b_0^2r_0^2=c_0$. We ensure this due to the
existence of a square root in the basic field $\BF$:
$r_0=\frac{\sqrt c_0}{b_0}$. The desired $r$ is then for example
$\frac{\sqrt c_0}{ b_0} t^{\frac{1}{2} {\rm deg\,}c-{\rm deg\,}b}$.

III.

In the case $\deg c = 2 \deg b$ we are again trying to find the
change of variables of the same type like in II, such that $\deg
(b^2 r^2+b^2 r+c) < \deg c$. Now $\deg r$ have to be zero and
cancellation of the highest terms is possible if and only if the
equation  $ b_0^2 r^2 + b_0^2 r +c_0=0$ solvable in $\BF$. If so,
after corresponding change of variables we get an equation with the
free term $c$ of the smaller degree. Otherwise, let us show that we
are in the imaginary case, i.e. there are no solutions of the
equation $\DD^2+b \DD+ c =0$ in $\BF((x))$. If we suppose that such
a solution does exist then $\deg \DD =\deg b$ and $ \DD_0^2  + b_0
\DD_0 +c_0=0$ ($\DD_0 \in \BF$ is a coefficient near the highest
term of $\DD$). But this means that the equation $ b_0^2 t^2 + b_0^2
 t +c_0=0$ is also solvable: $t=\frac{\DD_0}{b_0}$. Thus we are in
the imaginary case here if the equation is unsolvable in $\BF$.

IV.

Let $\DD \in K \backslash \BF [x]$ be the root of (\ref{t}):
$\DD^2+b\DD+c=0$.

If ${\rm deg\,} \DD^2={\rm deg\,}c$, then ${\rm deg\,} b\DD$ is greater.

We construct now the root in case ${\rm deg\,}b\DD = {\rm deg\,}c$,
i.e. ${\rm deg\,} \DD = {\rm deg\,}c-{\rm deg\,}b$. Denote $k={\rm
deg\,}\,\DD$, $m={\rm deg\,}b$, $\DD=a_k t^k+...$.  From the
equation $ \DD^2= b\DD+ c$ \,  we have $ a_k t^{2k} +...=a_k b_m
t^{m+k} +...+ c_{m+k} t^{m+k} +... $.  Since $2k<m+k$, for
cancellation it is necessary that $a_k b_m=c_{m+k}$, i.e.
$a_k=\frac{b_m}{c_{m+k}}$. Denote by $\widetilde\DD=\DD - a_k t^k$.
Then $\widetilde\DD$ satisfies the equation ${\widetilde\DD}^2 = b
\widetilde\DD+\widetilde c$, where $ \widetilde c = a_k^2 t^{2k} + b
a_k t^k + c$.  It is easy to see that ${\rm deg\,} \widetilde c <
{\rm deg\,} c $.  Hence we have to find the root of the equation
satisfying the condition in III, and this root will have the degree
smaller than $\BF$. By such an inductive procedure we obtain a
desired root as a power series.

In case ${\rm deg\,} \DD^2 = {\rm deg\,}b\DD$ we get a conjugate
root $b+\DD$.

\square

Consideration of the case II leads us necessarily  to the case I,
III or IV. Cases I and IV are imaginary and real respectively.
Consideration of the case III leads us either to the case I or IV,
or we remain in the imaginary case.

\section{Imaginary case}

Now we will give a solution in the imaginary case. According to
proposition \ref{Theorem 3.1.} we can assume that either 1). $ {\rm
deg\,}  c>2 {\rm deg\,} b$ and $ {\rm deg\,}  c$ is odd or 2). $\deg
c = 2 \deg b$ and the equation $ b_0^2 t^2 + b_0^2 t +c_0=0$
unsolvable in the field $\BF$.

Consider first the first case.

The main our tool here is a construction of the degree function on
$R$, which respects the multiplication.

\vspace{5mm}

\noindent {\bf Definition. } {\it Let define a function $ {\rm Deg}
:R \to {\mathbb Q}_+ \cup \{ -\infty \}$, as follows:
$${\rm Deg}(u+\Delta v)=\max
( {\rm deg\,}  u, {\rm deg\,}  v+\frac12  {\rm deg\,}  c),$$ where
{\rm deg\,} is the  usual degree function on polynomials. }

\vspace{5mm}

\begin{theorem}
\label{Theorem 4.1.} Let $f(t)=t^2+bt+c$ with $\rm {deg\,} \, c $ to
be odd and $R=\BF[x,t]/f(t)$, then for any $\alpha,\beta\in R$,
$$
{\rm Deg}(\alpha\beta)={\rm Deg}\alpha+{\rm Deg}\beta.
$$

\end{theorem}
{\bf Proof.} Let $\alpha=u'+\Delta v'$, $\beta=u+\Delta v$. Consider
four different possibilities for the degrees of $\alpha$ and
$\beta$:

1) ${\rm Deg} \alpha= deg\, u'$, ${\rm Deg} \beta={\rm \deg\,} \, u$ (i.e.
${\rm \deg\,} \, u'>\frac12{\rm \deg\,} \, c+{\rm \deg\,} \, v'$ and
${\rm \deg\,} \, u>\frac12{\rm \deg\,} \, c+{\rm \deg\,} \, v$);

2) ${\rm \Deg\,} \, \alpha=\frac12{\rm \deg\,} \, c+{\rm \deg\,} \,
v'$, ${\rm \Deg\,} \, \beta={\rm \deg\,} \, u$ (i.e., ${\rm \deg\,}
\, u'< \frac12{\rm \deg\,} \, c+{\rm \deg\,} \, v'$ and ${\rm
\deg\,} \, u> \frac12{\rm \deg\,} \, c+{\rm \deg\,} \, v$);

3) ${\rm \Deg\,} \, \alpha={\rm \deg\,} \, u'$, ${\rm \Deg\,} \,
\beta=\frac12{\rm \deg\,} \, c+{\rm \deg\,} \, v$ (i.e., ${\rm
\deg\,} \, u'>\frac12{\rm \deg\,} \, c+{\rm \deg\,} \, v'$ and ${\rm
\deg\,} \, u<\frac12{\rm \deg\,} \, c+{\rm \deg\,} \, v$);

4) ${\rm \Deg\,} \, \alpha=\frac{1}{2}{\rm \deg\,} \, c+{\rm \deg\,}
\, v'$, ${\rm \Deg\,} \, \beta=\ frac{1}{2} {\rm \deg\,} \, c+{\rm
\deg\,} \, v$ (i.e., ${\rm \deg\,} \, u'< \frac{1}{2} {\rm \deg\,}
\, c+{\rm \deg\,} \, v'$ and ${\rm \deg\,} \, u< \frac{1}{2}{\rm
\deg\,} \, c+{\rm \deg\,} \, v$).

Note that by definition
$$
{\rm \deg\,} \, \alpha \beta=\max\{{\rm \deg\,} \,(u'u+cv'v),\frac12{\rm \deg\,}
 \, c+
{\rm \deg\,} \,(u'v+v'u+bv'v)\}.
$$

In case~1 from the inequalities ${\rm \deg\,} \, u'>\frac12{\rm \deg\,} \, c+
{\rm \deg\,} \,
v'$,  ${\rm \deg\,} \, u>\frac12{\rm \deg\,} \, c+{\rm \deg\,} \, v$ and
${\rm \deg\,} \, c>2{\rm \deg\,} \, b$  we
have ${\rm \deg\,} \, u'u>{\rm \deg\,} \, cv'v$, ${\rm \deg\,} \, u'u>
\frac12{\rm \deg\,} \, c+{\rm \deg\,} \, u'v$,
${\rm \deg\,} \, u'u>\frac12{\rm \deg\,} \, c+{\rm \deg\,} \, v'u$, ${\rm \deg\,}
 \, u'u>\frac12{\rm \deg\,} \,
c+{\rm \deg\,} \, bv'v$. Hence ${\rm \deg\,} \, \alpha\beta={\rm \deg\,} \, u'u={\rm \deg\,}
\, u'+{\rm \deg\,} \,
u={\rm \deg\,} \,\alpha+{\rm \deg\,} \,\beta$.

In case~2 we similarly have ${\rm \deg\,} \,\alpha\beta=\frac12{\rm \deg\,} \,
c+{\rm \deg\,} \, u'v={\rm \deg\,} \, u'+\bigl(\frac12{\rm \deg\,} \, c+{\rm \deg\,} \,
v\bigr)={\rm \deg\,} \,\alpha+{\rm \deg\,} \,\beta$.

Case~3 is equivalent to case~2. One just has to replace $u'$ by $u$,
$u$ by $u'$, $v'$ by $v$ and $v$ by $v'$.

In case~4 we have ${\rm \deg\,} \,\alpha\beta={\rm \deg\,} \, cv'v=
(\frac12{\rm \deg\,} \, c+{\rm \deg\,} \, v')+(\frac12{\rm \deg\,} \, c+{\rm \deg\,}
\, v)={\rm \deg\,} \,\alpha+
{\rm \deg\,} \,\beta$.
\square

The existence of this degree function  allows us to solve the
equation $u^2+buv+cv^2=d$,  since it is equivalent to $(u+\Delta
v)(u+(b+\Delta) v)=d$ and therefore the degrees of $u$ and $v$
(which are non-negative) are bounded  by the degree of $d$.

Now we shell give a solutions in the second case. and show that
there are also a finite number of them.

Suppose that there exists a solution of (\ref{stard})  such that
$\deg u > \frac{1}{2} \deg d$ or $\deg v > \frac{1}{2} (\deg d- \deg
c)$. In this case degree of $d$ is not maximal, hence among $u^2,
buv, cv$ there are terms of the same degree. We show that if degrees
of two of them are coincide then the third has also the same degree.
It is easy calculations in three possible cases using $\deg c = 2
\deg b$. Hence the highest terms of $u^2, buv$ and $ cv$  are
cancelled and $u_0^2+b_0u_0v_0+c_0v_0^2=0$ holds for
$u_0,v_0,b_0,c_0 \in \BF$ -- coefficients near  the highest terms of
polynomials $u,v,b,c$. But it means that the equation $b_0^2
t^2+b_0^2 t + c_0=0$ is solvable: $t=u_0/v_0b_0$. Therefore we have
a bound for the degrees of $u,v$ also in this case.

We will use later on the following denotation:
$$r_{A,B}=\min\limits_{U\in \Gl(2): UAU^{ -1}=B}  \deg U.$$

In both variants of the imaginary case we get the following linear
estimation.

\begin{proposition}
\label{Proposition 4.1.} An estimation for the degree of elements of
the conjugating matrix in imaginary case is linear: $r_{A,B} \leq 2
\dd$, where as earlier $\dd$ is a maximum of $\deg (A)$ and $\deg
(B)$.
\end{proposition}
{\bf Proof.}
To obtain the estimation we have to take into account that before
we turn out to be in real or imaginary case we have to make a
change of variables of the type
$$\left\{\begin{array}{l}
u'=u+qv,\\
v'=v\end{array}\right.
$$
where $\deg q \leq \frac{1}{2} \deg c$.

Then  in imaginary case of  type I (proposition \ref{Theorem 3.1.})
we estimate degrees of $u'$ and $v'$ from the equality
$(u'+ \DD v')(u'+ (b+ \DD) v')=d$ using introduced above
degree function $\rm{Deg}$ on $R$. We get $\deg u' \leq \dd/2$,
$\deg v' \leq \dd/2$ and $\deg u \leq \dd$, $\deg v \leq \dd$.

In imaginary case of  type III (proposition \ref{Theorem 3.1.})
as was shown above we have bounds:
$\deg u \leq \frac{1}{2} \deg d$ or $\deg v \leq \frac{1}{2} (\deg d- \deg c)$.
Hence, also $\deg u \leq \dd$ and $\deg v \leq \dd$.

Now taking into account (\ref{pq})
 we have an estimation  for the degree of entries of conjugating
matrix: $r_{A,B} \leq 2 \dd$. \square

\section{Real case}\label{real}
\subsection{Units (equation $u^2+buv+cv^2=1$)}\label{runits}

We start in real case ($ {\rm deg\,}  c < 2 {\rm deg\,}  b$)
with the solution of
our equation with $d=1$:
\begin{equation}
u^2+buv+cv^2=1.
\label{star}
\end{equation}

If $(u,v)$ is a solution of our equation we  say also that $\oo=
u+\Delta v \in R$ is a solution. Denote by $U(R)$ the set of all
solutions $\oo \in R$ of the equation (\ref{star}), $U(R)$ becomes a
group with the multiplication by that of $R$.

\noindent {\bf Definition. } {\it We say that $p\in R$ is reduced if
$ {\rm deg\,}  p>0$ and $ {\rm deg\,}  p'<0$, where $p'$ is the
conjugate element (defined in section \ref{prelim})}.

\vspace{5mm}

\begin{theorem}
\label{Theorem 5.1.1.}
The set $U(R)$ of solutions of
(\ref{star})
is an infinite cyclic group.
The generator of $U(R)$ is an element with
 minimal positive degree.
Moreover $R^*=U(R) \times \BF^*$, where $R^*$ is the group of units
of $R$.

\end{theorem}
{\bf Proof}. Show first that $R^* =U(R) \times \BF^*$. The equation
(\ref{star}) means that $N(\oo)=1$, hence lemma \ref{Lemma 2.1.} c).
and d). implies  that $U(R)$ is a subgroup of $R^*.$ Let $\oo \in
R^*$. According to lemma \ref{Lemma 2.5.} $N(\oo) \in \BF^*$. Since
$\BF$ is finite field of characteristic $2$ there exists a unique
$\a \in \BF^*$ such that $N(\oo)=\a^2$. Therefore $N(\oo / \a)=1$,
i.e. $\oo / \a \in U(R)$ and $\a$ is uniquely determined.

Now we prove that $U(R)$ is an infinite cyclic group and its
generator  is an element with minimal positive degree.

\begin{lemma}
\label{Lemma 5.1.1.}
If $\ee \in R^*$  and $\vert \ee \vert =1$,
than $\ee$  is a nonzero constant.
\end{lemma}
{\bf Proof}. According to lemma \ref{Lemma 2.5.} $\vert N(\ee) \vert
=1$. Hence $\vert \ee' \vert= 1$ follows from $\vert N(\ee) \vert
=\vert \ee \ee' \vert= \vert \ee \vert \vert \ee' \vert = 1$.
Comparing  the corresponding power series we can see that $\vert \ee
\vert = 1$ and $\vert \ee' \vert= 1$ together imply that $bv \in \BF
$. Hence there are three possibilities: $b=0$; $v=0$; or  $b, v \in
\BF^*$. The case $b=0$ was considered in section \ref{prelim}, $v=0$
means that $\ee$ is a polynomial but it was a unit, so it is
actually a constant. From $b \in \BF^*$ and $v \in \BF^*$ it follows
that $c=0$ (since ${\rm deg\,} c < 2{\rm deg\,} b$) and we are in
the case when the equation is factorizable over $\BF[x]$, which
again was considered in section \ref{prelim}. \square

\begin{lemma}
\label{Lemma 5.1.2.} If $\ee_1$ and $\ee_2$  are units and $\vert
\ee_1 \vert =\vert \ee_2 \vert$, than $\ee_1$ and $\ee_2$ coincide
up to a constant: $\ee_1= \a \ee_2, \a \in \BF^*$.
\end{lemma}
{\bf Proof}. Obviously $\frac{\ee_1}{\ee_2}$ is also a unit and
$\vert \frac{\ee_1}{\ee_2} \vert = \frac{\vert \ee_1 \vert}{ \vert
\ee_2 \vert} =1$, hence  by lemma \ref{Lemma 5.1.1.}
$\frac{\ee_1}{\ee_2} \in \BF^*$. \square

Let $\ee_0$ be the unit with minimal valuation $\vert \ee \vert > 1$
(with minimal positive degree).

\begin{lemma}
\label{Lemma 5.1.3.} Any unit $\ee \in R^*$ has the form $\ee = \a
\ee_0^n, \a \in \BF^*$
\end{lemma}
{\bf Proof}.
Suppose that it is not true. There exists  $n \in {\Bbb N}$ such that
$ {\vert \ee_0 \vert}^n < \vert \ee \vert
< {\vert \ee_0 \vert}^{n+1}$.
The equality is impossible, because if
$\vert \ee \vert= {\vert \ee_0 \vert}^{n}$
than by lemma \ref{Lemma 5.1.2.} $\ee = \a \ee_0^n$.
We then multiply previous inequalities by
${\vert \ee_0 \vert}^{-n}$ and get a
contradiction with minimality
of $\vert \ee_0 \vert$:
$ 1< \vert \ee_0^{-n} \ee \vert < \vert \ee_0 \vert$.
\square

By this the proof of the theorem is completed.

\square

\vspace{5mm}

We find the generator of $U(R)$ in two steps.
First, we
construct some nontrivial element of $U(R)$.

Let us denote by  $[A_0;A_1,A_2,...] $ where $ (A_i \in \BF[x])$ the
continued fraction expansion $A_0+\frac{1}{A_1+\frac{1}{A_2+...}}$.
We shall say that this expansion is {\it purely} periodical if the
periodicity of the sequence $A_0;A_1,A_2,...$ starts from $A_0$.

\vspace{5mm}

\begin{theorem}
\label{Theorem 5.1.2.} Let $\p\in R$ be a reduced root of (\ref{t}).
Then the continued fraction expansion $\p=[A_0;A_1,A_2,...] \, (A_i
\in \BF[x])$ is purely periodical with a period $T \leq q^{2m}$,
where $q=|\BF|, m={\rm deg\,}b.$
\end{theorem}
{\bf Proof}. One can present the series $\p_n=[A_n;A_{n+1},...]$
which appears in the process of construction of a continued
fraction, as obtained by operations $\phi_1: \p \rightarrow u+ \p$
(cutting a polynomial part of the series) and $\phi_2: \p
\rightarrow 1/ \p$ (taking an inverse).

It is easy to see that $\phi_1$ and $\phi_2$ act on the set ${\cal
U} = \{ $solutions of the equations $\wt ax^2 + bx + \wt c=0 \ \
\vert {\rm deg\,}  \ \  \wt c < {\rm deg\,} b, \ \ {\rm deg\,} \wt a
< {\rm deg\,} b\}$.

Indeed, let $x \in {\cal U}$, $y=\phi_1(x)=x+u$, where $u \in
\BF[x], \, {\rm deg\,} y < 0 $. Since $\wt ax^2 + bx + \wt c=0$, we
have that $\wt a y^2 + by + c'=0$, where $c'= au^2 + bu + \wt c$.
From the latter equation $c'= \wt a y^2 + by$, and ${\rm deg\,} y <
0 $. Hence for the degree of $c'$ we have ${\rm deg\,} c' < {\rm
deg\,}(\wt a y+b) \leq {\rm deg\,} b$, therefore $y \in {\cal U}$.

Let now $x \in {\cal U}$ and $y=\phi_2(x)=1/x$. Since $\wt ax^2 + bx
+ \wt c=0$, we have that $\wt cy^2 + by + \wt a=0$, therefore $y \in
{\cal U}$. Thus  $\phi_1(x)$ and $\phi_2(x)$  acts on ${\cal U}$.

Hence the number of steps to obtain the same $\p$ is less then
$|{\cal U}|=2q^{2m}$, where $m= {\rm deg\,}b, q$
--- number of
elements of the field.

Put now ${\cal U}_+ = \{x \in {\cal U}, {\rm deg\,} x \geq 0\}$.
Note that $\p_{n+1}= \phi_2 \phi_1(\p_n)$. Purely periodicity
follows from the fact that $\phi_2 \phi_1$ is a permutation of
${\cal U}_+$. Indeed, $\phi_2 \phi_1(\p_n)$ acts on ${\cal U}_+$ and
it can be easily checked that it is an injection.

The estimation $T<q^{2m}$ follows from the equality $|{\cal
U_+}|=q^{2m}.$

\square

 \vspace{5mm}

For  $n=T$ we have $\p=\p_T$, and
$$
\p=\frac{P_n\p+P_{n-1}}{Q_n \p + Q_{n-1}},
$$
where

\begin{equation}
P_{n+1}=P_n A_n+P_{n-1},  P_0=1,  P_1=A_0
\label{recP}
\end{equation}

\begin{equation}
 Q_{n+1}=Q_n A_n+Q_{n-1}, Q_0=0,  Q_1=1
\label{recQ}
\end{equation}

\vspace{5mm}

This means that $\p$ satisfies the quadratic equation $Q_n
\p^2+(P_n+Q_{n-1})\p+P_{n-1}=0$. Since $\p$ satisfies also the
equation $\p^2+b\p+c=0$ and the latter equation does not have
solutions in rational functions (this case was considered separately
in the section \ref{prelim}), these two equations are proportional.
Denote the coefficient of proportionality by $V$. Then $Q_n=V$,
$P_{n-1}=cV$, $P_n+Q_{n-1}=bV$. Denote $P_n=U$. From the known
equation $P_nQ_{n-1}+Q_nP_{n-1}=1$ we obtain that $\epsilon=U+\Delta
V=P_T+\Delta Q_T$ is a solution of the equation (\ref{star}).

\vspace{5mm}

\begin{lemma}
\label{Lemma 5.1.4.} When we live in the real case ($ {\rm deg\,} c<
2{\rm deg\,}  b$), there exists an invertible linear change of
variables which turns the equation $u^2+buv+cv^2=d$ into
$u^2+buv+\widetilde cv^2=d$ with ${\rm deg\,} \widetilde c<{\rm
deg\,} b$.
\end{lemma}
{\bf Proof}.
We have $u^2+buv+cv^2=d$,
$ {\rm deg\,} b \leq {\rm deg\,} c< 2{\rm deg\,}  b$. Let us divide $c$
by $b$: $c=bq+r, {\rm deg\,} r < {\rm deg\,} b$ and consider
the change of variables:
$$\left\{\begin{array}{l}
u'=u+qv,\\
v'=v\end{array}\right.
$$
New equation is: $u^2+buv+(q^2+r)v^2=d$. Denote $\widetilde
c=q^2+r$. Note that ${\rm deg\,} \widetilde c < {\rm deg\,} c$.
Indeed, ${\rm deg\,} r < {\rm deg\,} b \leq {\rm deg\,} c$ and ${\rm
deg\,} q={\rm deg\,}c-{\rm deg\,}b$, hence ${\rm deg\,}q^2=2{\rm
deg\,}c-2{\rm deg\,}b<{\rm deg\,}c$, this means that we can keep
making changes of variables of such a type until we get ${\rm
deg\,}b>{\rm deg\,}\widetilde c$. \square

Note that the composition of changes of variables of the type
$u'=Qu, v'=v$, with $deg\, Q \leq \alpha$ has the same form.

\vspace{5mm}

In proposition \ref{Theorem 3.1.} III  we proved that the case $
{\rm deg\,}  c<2 {\rm deg\,}  b$ is real by the construction of the
root $\Delta$ of the equation (\ref{t}) as a power series. It
follows from this construction that in the case $ {\rm deg\,} c<{\rm
deg\,} b$ one of the roots of our equation is reduced. Hence we have
obtained the following lemma.

\vspace{5mm}

\begin{lemma}
\label{Proposition 5.1.1.} If $ {\rm deg\,}  c< {\rm deg\,}  b$ then
the element $\Delta+b$ is reduced and $\epsilon=P_T+\Delta Q_T$ is a
nontrivial element of the group $U(R)$, where $T$ is the period of
the continued fraction expansion of $\Delta+b$.
\end{lemma}

\begin{lemma}
\label{Proposition 5.1.2.} The estimation for the degree of this
element $\epsilon$ of the group $U(R)$ is the following: ${\rm
deg\,} \epsilon \leq  \dd q^{2 \dd}$.
\end{lemma}
{\bf Proof}. We have to estimate first $\deg u' = \deg P_T$ and
$\deg v' = \deg Q_T$. Here $u'$ and $v'$ are the same as at the
proof of proposition \ref{Proposition 4.1.}.
 Recall that $\p=[A_0;A_1,A_2,...]$ is a continued fraction
expansion of the reduced root of (\ref{t}). Note that $\deg A_n \leq
\deg b$. Indeed, it is a positive part of an element
$\p_n=[A_n;A_{n+1},...] \in \cal U$ ($\cal U$ is the set constructed
above in the proof of \ref{Theorem 5.1.2.}) and $\deg \p_n=\deg b -
\deg \widetilde c \leq \deg b $. From the recurrent formulas
(\ref{recP}) and (\ref{recQ}) for $P_n$ and $Q_n$ it follows that
$\deg u' \leq \deg b \, T \leq \dd q^{2 \dd}$ and $\deg v' \leq \deg
b \, (T-1) \leq \dd (q^{2 \dd}-1)$. Then we get estimations for
$\deg u$, $\deg v$ and $\deg \ee=\deg (u + \DD v) \leq \dd q^{2
\dd}$. \square

\vspace{5mm}

Now we have to construct the generator of U(R), from a nontrivial
element of U(R), we just have found using the continued fraction
expansion.

\vspace{5mm}

\begin{lemma}
\label{Theorem 5.1.4.}
Let $\epsilon_0=x_0+\Delta y_0$ be any generator
of the group $U(R)$. Then $y_0=g.c.d.(Y)$, where
$Y=\{y:x+\Delta y\in U(R)\}$.
\end{lemma}
{\bf Proof}.
Let $\oo \in U(R), \oo=x+\DD y, \oo^n=x_n+\DD y_n$.
It is enough to check the following recursive formula:
$y_{n+2}=y_n+by y_{n+1}$.
\square

\vspace{5mm}

Hence we can just consider all divisors $y_i$ of $y$ where
$\epsilon=x+\Delta y$ is an element of $U(R)$, we have constructed.
Then find $x_i$, such that $x_i^2+bx_iy_i+cy_i^2=1$, for those $y_i$
for which it is possible. From the constructed in this way finite
set of $\epsilon_i\in U(R)$ we select those with minimal positive
degree. This is the desired $\epsilon_0$.

We can summarize the results of this section in the
following

\begin{proposition}
\label{Proposition 5.1.5.} The group $U(R)$ is an infinite cyclic
group and there exists an algorithm for constructing of its
generator.
\end{proposition}

\vspace{5mm}

\subsection{General case $d\neq1$}

Now we consider a general equation (\ref{stard}):
$u^2+buv+cv^2=d, d \neq 1$.

Let $\epsilon_0$ be the generator of the group $U(R)$ with a
positive degree. Denote $k= {\rm deg\,} \epsilon_0$.

\vspace{5mm}

\begin{lemma}
\label{Theorem 5.2.1.} The set of all solutions of (\ref{stard}) has
the form $\{\omega\epsilon_0^l:l\in \mathbb Z, \, \omega$ is a
solution of (\ref{stard}) with $ {\rm deg\,} \omega=0,\dots,k-1\}$.
\end{lemma}
{\bf Proof}. If $\oo$ is an arbitrary solution of (\ref{stard}):
$N(u,v)=u^2+buv+cv^2=d$, then any $\oo \ee_0^l, l \in {\Bbb Z}$ is
also a solution (lemma \ref{Lemma 2.1.} c).  Hence we can rewrite
the set of solutions of (\ref{stard}) in the following way: $\{\oo
\ee_0^l: \, l \in {\mathbb Z}, \oo$ is a solution of (\ref{stard})
with ${\rm deg\,} \oo =0,1,...,k-1 \}$.  \square

\vspace{5mm}

\begin{theorem}
\label{Theorem 5.2.2.}
Let $\omega=u+\Delta v$ be a solution of
(\ref{stard}) with  ${\rm deg\,} \omega=0,\dots,k-1$. Then $
{\rm deg\,}  v\leq \max\{ {\rm deg\,}  d,k\}- {\rm deg\,}  b$.
\end{theorem}
{\bf Proof}. As we said at the beginning of the section \ref{runits}
\, $\oo=u+\DD v$ is a solution of  (\ref{stard}) means that $(u,v)$
is a solution of (\ref{stard}), i.e. $N(\oo)=\oo \oo'=d$. On the
other hand $\oo \oo'= \oo^2 +b \oo v$. Hence the solution $\oo$
satisfies the equation
\begin{equation}
\oo^2+(bv)\oo+d=0.
\label{om}
\end{equation}

Let us consider two cases:

Case 1. ${\rm deg\,} bv \leq {\rm deg\,} d$,

For this case ${\rm deg\,} v \leq {\rm deg\,}d-{\rm deg\,}b.$ \\

Case 2. ${\rm deg\,} bv > {\rm deg\,} d$,

Here a priori there exist
three possibilities for the degrees of terms in the equation.

a). ${\rm deg\,} \oo^2 = {\rm deg\,}d$.
This is impossible because it implies
${\rm deg\,} bv > {\rm deg\,} \oo^2={\rm deg\,} d$
and the highest term  $(bv)\oo$ can not be cancelled.

b). ${\rm deg\,} \oo^2 = {\rm deg\,} bv\oo$.
It means that ${\rm deg\,} \oo = {\rm deg\,} bv$. Since we are
interested in the solutions $\oo$ with ${\rm deg\,} \oo \leq k-1$,
we have ${\rm deg\,} bv \leq k-1$,
and ${\rm deg\,} v \leq k-1- {\rm deg\,} b$. This proves the theorem in this case.

c). ${\rm deg\,} bv\oo = {\rm deg\,} d$.
Hence ${\rm deg\,} \oo < 0$. Which is incompatible with the hypothesis.

We can conclude that for any solutions $\oo$ of (\ref{stard})
we have
${\rm deg\,} v \leq {\rm max} ({\rm deg\,} d -
{\rm deg\,}b, k-1-{\rm deg\,}b).$

\square

According to the theorem \ref{Theorem 5.2.2.}  we can find all
solutions of (\ref{stard})  by a finite procedure.

\vspace{5mm}

The last step is to choose from the set of solutions
$\omega\ee_0^l=u_l+\Delta v_l$ those for which $b_{21}$ is a divisor
of $u_l$, and  $a_{21}$ is a divisor of both $(a_{11}+b_{11}) u_l +
b_{12} v_l$ and $b_{12} u_l+(a_{11}+b_{22}) v_l $. It is possible to
describe all such solutions due to the following fact.

\vspace{5mm}

\begin{lemma}
\label{Theorem 5.2.3.}
Let $\omega\ee_0^l=u_l+\Delta v_l$,
$P$ be a polynomial and $r_l$,
${\rm deg\,} r_l < {\rm deg\,}P$, be the sequence of
residues of $P_1 u_l+P_2 v_l$ for some polynomials $P_1,P_2$
with respect to $P$. Then this sequence is
periodical: $r_{l+T_0}=r_l$  for some period $T_0$.

The estimation for the period is the following: $T_0 \leq q^{\deg
P}$.
\end{lemma}
{\bf Proof.} It is enough to show the periodicity of residues of
$x_n$ and $y_n$, where $\ee_0=x+\DD y$,    $ \ee_0^n=x_n+\DD y_n$.
Let $\res (Q,P)$ denote the residue of $Q$ with respect to $P$. It
is clear that $\res (Q(x_n,y_n),P)= \res Q(res(x_n,P), \res
(y_n,P))$. We will show that $\res (x_n,P)$ and $\res (y_n,P)$ are
periodical with the period $T_0$. Then periodicity of $x_n$ and
$y_n$ will follow, because
 $\oo \ee_0^n=u_0x_n+cv_0y_n + \DD (u_0 y_n +v_0 x_n+bv_0 y_n)$,
where  $\oo = u_0+ \DD v_0$.
Let $r_n=\res (y_n, P)$,  $ s_n=\res (x_n, P)$.
Just from $(x_n+\DD y_n)(x+\DD y)=x_{n+1}+\DD y_{n+1}$,
we have the following recurrent formulas:
$$ x_{n+1}=x_n x + c y_n y, \, \\
y_{n+1}=x_n y + y_n x + b y_n y.$$
We shell consider the sequence of pairs of the residues: $(r_n,s_n)$.
It is recurrent of length one since
$$ s_{n+1}=\res (s_n x+c r_n y), \, \\
r_{n+1}=\res (s_n y+ r_n x +b r_n y).$$
This sequence belongs to the finite set $M^2$
of pairs of polynomials of degree $ < \deg P$,
hence it is periodical with a period $T_0 \leq q^ {\deg P}$.
\square

As a corollary we see that the estimation for the period $T_0$ in
our situation is the following: $T_0 \leq q^{2 \dd}$, where $\dd$ is
a maximum of degrees of entries of initial matrices.

Summarizing above statements of
lemma \ref{Theorem 5.2.1.}, theorem \ref{Theorem 5.2.2.}
and lemma \ref{Theorem 5.2.3.} we end up with a construction of
the set $\wt{ \sm}_{A,B}$ which describes the set of all
conjugating matrices for the pair  $A$, $B$.
$$
\begin{array}{l}
\wt{ \sm}_{A,B}=\{ \oo \ee_0^{l+nT_0} \vert  {\rm deg\,} \oo =
0,...,k-1, \ \
l=0, \dots T_0-1, \\
b_{21} \Big| u_l; \ \  a_{21} \Big|
(a_{11}+b_{11}) u_l + b_{12} v_l;  \ \ {\rm and} \
\ a_{21} \Big| b_{12} u_l+(a_{11}+b_{22}) v_l \}.
\end{array}
$$
For any element $\oo \in \wt{ \sm}_{A,B}$ one can obtain a
conjugating matrix
$$
U=\left( \begin {array} {cc}
                u&p\\
                v&q
        \end {array} \right),
$$
where $\oo=u+ \DD v$ and $p,q $ found from (\ref{pq}).

Moreover from $\wt{ \sm}_{A,B}$ we can also choose a finite subset
which characterize the conjugacy of matrices $A$ and $B$.

{\bf Corollary.} {\it Put
$$
\begin{array}{l}
 \sm_{A,B}=\{ \oo \ee_0^l \vert  {\rm deg\,} \oo =
0,...,k-1, \ \
l=0, \dots T_0-1, \\
b_{21} \Big| u_l; \ \  a_{21} \Big|
(a_{11}+b_{11}) u_l + b_{12} v_l;  \ \ {\rm and} \
\ a_{21} \Big| b_{12} u_l+(a_{11}+b_{22}) v_l \}.
\end{array}
$$
then $\sm_{A,B} \neq \emptyset \Longleftrightarrow A$
is conjugate with $B$.}

\begin{proposition}
\label{Proposition 5.2.2.} The estimation for the degree $r_{A,B}$
of the entries of conjugating matrix in real case is the following:
$r_{A,B} \leq
2 \dd q^{6 \dd}$.
\end{proposition}
{\bf Proof.}
We will estimate first degree of the solution $\widetilde{\oo}= \oo \ee_0^l$
where $\deg \oo \leq k-1$, $l \leq T_0-1$.

Let $\oo \ee_0^l=u_l+\DD v_l, T_0={\rm l.c.m} (T_1,T_2,T_3)$, where
$T_1,T_2$ and $T_3$ are periods of the residues (as they defined in
lemma \ref{Theorem 5.2.3.}) of sequences $u_l$, $(a_{11}+b_{11}) u_l
+ b_{12} v_l$ and $b_{12} u_l+(a_{11}+b_{22}) v_l $ relative to
$b_{12}, a_{21}$ and $ a_{21}$ respectively. Using proposition
\ref{Theorem 5.2.3.}
 we can estimate $l.c.m(T_2,T_3)$
by $ q^{2 \deg a_{21}}$, $T_1$ by $ q^{2 \deg b_{12}}$, hence $T_0
\leq q^{4 \dd}$. Then $\deg \widetilde{\oo}= \deg \oo \ee_0^l \leq
\dd q^{6 \dd}-1$. Using the equation (\ref{om}) we can obtain: $\deg
u \leq \dd q^{6 \dd}+ \dd -1$ and $\deg v \leq \dd q^{6 \dd}+ \dd$.
(Note that in real case it is impossible that $\dd =0$.) \square

Combining together all  estimations obtained above separately in the
following cases: when one of the matrices is diagonal,
 $b=0, \DD \in \BF(x)$, real case and imaginary case, we
get the estimation for the degree of conjugating matrix in case of
$char=2: \quad r_{A,B} \leq \dd(q^{6\dd} +2)$.

\vspace{8mm}

\section{ Note on the case of characteristic $\neq 2$ }

\vspace{5mm}

The case of positive characteristic $\neq 2$ was considered in
Artin's paper \cite{EA}. We will obtain here only an estimations
which comes not always from the procedure described in \cite{EA}. We
also have to make precise here the reduction of the conjugacy
problem to a  quadratic equation in this case.

Let $char \BF \neq 2$. The conjugating condition $UAU^{-1} = B$  for
$U \in \Gl(2,\BF[x])$ is equivalent to the system
$$
                uq-pv \in {\BF}^*
$$
$$
                UA=BU
$$
where $ U=\left( \begin {array} {cc}
                u&p\\
                v&q
        \end {array} \right)$,
${\BF}^*=\BF \backslash \{0\}$ is the multiplicative group of the
field $\BF$.

To check the solvability of the system for an arbitrary coefficient
$\alpha \in {\BF}^*$ we need to check it only in two cases:
 $\alpha =1$ and $\alpha$ is an element of ${\BF}^*$ which is not a square.
(Other solutions can be obtained from them because $\alpha / \beta$
is a square if $\alpha$ and $\beta$ are not.) We will obtain then a
quadratic equation of the type $u^2+b uv + v^2=d$ in the similar way
as earlier, and since $char \neq 2$ we are able to reduce it to the
Pell's equation: $u^2-cv^2=d$.

So we can use here the conventional notions of real and imaginary
case and usual rule to distinguish them:

\begin{proposition}\label{Proposition 6.2.2.}\ \break

I. If $ \deg \, c$ is odd or  $ \deg \, c$ is even but the highest
term of polynomial $c: c_n$ is not a square, then we are in the
imaginary case.

\vspace{2 mm}

II. If $\deg \, c$ is even and $c_n$ is a square then we are in the
real case.
\end{proposition}

To obtain the estimations in these cases we will treat them
separately.

In the imaginary case the highest terms of $c$ is not a square of
 a monomial (with coefficient).
If we suppose that there exists a solution with $\deg \, u^2 > \deg
\, d$ or $\deg \,c v^2 > \deg \, d$ then the highest terms of $u^2$
and $c v^2$ have to be cancelled, but it would  mean that highest
term of $c$ is a square of a monomial (with a coefficient). Hence
for any solution $\deg u < \frac{1}{2} \deg d$ and $\deg v <
\frac{1}{2} (\deg d - \deg c)$. And we have got the following
estimation.

\begin{theorem}
\label{Proposition 6.2.3.} In the imaginary case always there exist
only a finite number of solutions, and estimation for the degree of
entries of conjugating matrix is linear:  $r_{A,B} \leq 2 \dd$,
where $\dd$ is a maximum of $\deg (A)$ and $\deg(B)$.
\end{theorem}

Now consider the real case. In Artin's paper it was proved (in case
of $char \neq 2$)  that there exists a reduced root $\p$ (the root
with $\deg \p
> 0$ and $ \deg \p'<0$ for the conventional definition of $\p'$)
of the equation $t^2+bt+c=0$ and the continued fraction expansion of
this reduced root is purely periodical. Also he proved that if
continued fraction expansion is purely periodical then some unit can
be constructed. He does not give an estimation of the period. We get
it here now by methods similar to those we have used above, with
only few essential changes. Namely, we present the process of
continued fraction construction via an actions on the finite set.
The set $\cal U$ we have to take  to deal with the case $char \neq
2$ is different from one appeared in the proof of the Theorem
\ref{Theorem 5.2.1.}. Estimation for the size of this set gives an
estimation of the period.

\begin{lemma}
\label{Lemma 6.2.1.} Let $\p \in R$ be a reduced root of the
equation $at^2+bt+c=0$ with the condition: ${\rm deg\,} a \leq \ggg,
{\rm deg\,} b \leq \ggg$ and ${\rm deg\,} c \leq \ggg$, where $\ggg$
is as earlier a maximum of degrees of elements of initial matrices
$A$ and $B$.  Then the continued fraction expansion $\p
=[A_0;A_1,A_2,...] \, (A_i \in \BF[x])$ is periodical with a the
period $T \leq q^{3\ggg}$.
\end{lemma}
{\bf Proof}. One can present the series $\p_n=[A_n;A_{n+1},...]$
which appears in the process of construction of a continued fraction
expansion, as obtained by the operations $\phi_1: \p \rightarrow \p
- p $ (cutting a polynomial part of the series) and $\phi_2: \p
\rightarrow 1 / \p$ (taking an inverse).

We show that $\phi_1$ and $\phi_2$ act on the set ${\cal U}_r$  of
reduced roots of the following equations: ${\cal U}_r = \{ $reduced
roots of the equations $ at^2 +  bt +  c=0 \ \ \vert \ \ {\rm deg\,}
a \leq \ggg, \ \ {\rm deg\,}  b  \leq \ggg, \ \ {\rm deg\,} c \leq
\ggg \}$. But the root can be reduced only if $\deg b - \deg a < 0$
and $\deg c - \deg b < 0$    hence in fact ${\cal U}_r = \{ $reduced
roots  of the equations $ at^2 +  bt +  c=0 \ \ \vert \ \ {\rm
deg\,}    a < \ggg, \ \ {\rm deg\,}   b  \leq \ggg, \ \ {\rm deg\,}
 c <  \ggg \}$. It is known (\cite{EA}) that if in the process of
continued fraction expansion we get reduced root $\p_n$, then all
$\p_{n+k}$ for $ k\in \mathbb Z$ will be also reduced roots.

Hence we have to show only that if we take $\p \in {\cal U}_r$, then
$\phi_1(\p) $ and $\phi_2(\p) \in {\cal U}$, here ${\cal U} = \{
$solutions of the equations $ at^2 +  bt +  c=0 \ \ \vert \ \ {\rm
deg\,}  a \leq \ggg, \ \ {\rm deg\,}  b  \leq \ggg, \ \ {\rm deg\,}
c \leq  \ggg \}$.

It is obvious for $\phi_2$. Let $ \p \in {\cal U}$ and
$y=\phi_2(\p)=1/\p$. Since $ a\p^2 +  b\p +  c=0$, we have $ cy^2 +
by +  a=0$ and  $y \in {\cal U}$.

To prove the same property for $\phi_1$ find first how the equation
changes. If $\phi_1(\p) = \p - p =y$ then for $y$ we have: $a y^2 +
(2ap+b)y+(ap^2+bp+c)=0$.
A priori there are three possibilities: 1). ${\rm deg\,} a{\p}^2 =
{\rm deg\,} c > {\rm deg\,} b \p$; 2). ${\rm deg\,} a{\p}^2 = {\rm
deg\,} b \p > {\rm deg\,} c$; 3). ${\rm deg\,} b \p = {\rm deg\,} c
> {\rm deg\,} a{\p}^2$, but only the second one could actually exist. In this
case ${\rm deg\,}(2ap+b) = {\rm deg\,}b $ and ${\rm deg\,}
(ap^2+bp+c) \leq \ggg$. The latter follows from $ap^2+bp+c=a( \p
-y)^2+b( \p -y)+c=a{\p}^2+b{\p} + c + 2a \p y - by + a y^2$
and ${\rm deg\,} y < 0$. Hence we have the equation on $y$ of the
same type and $y \in {\cal U}$.

Notice that the highest terms of $b$ (middle coefficient) are the
same for all elements of ${\cal U}_r$. Hence we can estimate the
number of elements of ${\cal U}_r$ as follows: $|{\cal U}_r| \leq
q^{3 \ggg}$. It is an estimation for the period of continued
fraction expansion. \square

Now we shall obtain the final estimation of $r_{A,B}$ in the case of
$char \neq 2$ based on the estimation of the period.

\begin{theorem}
\label{Proposition 6.2.4.} In the real case the estimation for the
degree of entries of the conjugating matrix is the following:
$r_{A,B} \leq
(q+1)\dd q^{7 \dd}$, where $\dd$ is a maximum of  $\deg (A)$ and
$\deg(B)$.
\end{theorem}
{\bf Proof. } First we consider the root $\DD$ of the equation
$t^2=c$. We can construct it as a power series by usual recursive
procedure. It is not a reduced root, since here $\deg \DD= \deg
\DD'$.
 Note that it is different from the case of $char =2$
where one of two roots of the initial equation had to be reduced.
But in process of the construction of the continues fraction
expansion of $\DD=[A_0;A_1,...]$ in some step
$\DD_n=[A_n;A_{n+1},...]$  the reduced root have to appear
(\cite{EA}). For this root $\p=\DD_{n_0}$ we will have
$\p=\frac{P_n\p+P_{n-1}}{Q_n \p + Q_{n-1}}$. Let $\p$ satisfy the
equation: $A \p^2+ B \p + C=0$. Then from the proportionality of two
quadratic equations (which can not have a rational solutions) for
the reduced root, we get a formulas for nontrivial unit $X+\DD Y$:
$P_n+Q_{n-1}=2X, \, P_n-Q_{n-1}=2Yb$, where $n=T$, $T$ is a period
of the continued fraction expansion of $\p$. Hence
 we can estimate
degrees of $X$ and $Y$ as $\deg X \leq \deg B \, T$, $\deg Y \leq
\deg B \, T$. Let us convince that $\deg B \leq \dd$ and $T \leq
q^{3\dd+1}$, for $\dd$ being as earlier maximal degree of entries of
given matrices. It will follow from lemma \ref{Lemma 6.2.1.}. But we
have to note the following: if we start the process of continued
fraction expansion  with $\DD$, the root of $t^2=c$, $\deg c \leq 2
\dd$, then already at the first step, we get an equation
$\DD_1^2-2p\DD_1-c+p^2=0$ with $\deg(c-p^2) \leq \dd, \deg(2p) \leq
\dd$. At the next steps degrees of the coefficients of the equations
can not become bigger any more.

Hence we can estimate $\deg X$ and $\deg Y$ by $\dd q^{3\dd+1}$.
If we take into account the change of variables,
we get $ \dd q^{3 \dd+1} +\dd$.

Then the estimation  for $k=\deg\, \ee_0, \ee_0 = X+ \DD Y$
will be: $k \leq  \dd q^{3 \dd+1} + 2\dd$.

Let $\oo \ee_0^l=u_l+\DD v_l, T_0={\rm l.c.m} (T_1,T_2,T_3)$, where
$T_1, T_2$ and $T_3$ are periods of the residues (as they defined in
lemma \ref{Theorem 5.2.3.}) of sequences $u_l$, $(a_{11}+b_{11}) u_l
+ b_{12} v_l$ and $b_{12} u_l+(a_{11}+b_{22}) v_l $ relative to
$b_{12}, a_{21}$, and $a_{21}$ respectively. By the estimation from
the proposition \ref{Theorem 5.2.3.} we have $T_0 \leq q^ {4 \dd}$.

Consider an arbitrary root of $u^2 + c v^2=d: \, \widetilde \oo= \oo
\ee_0^l= u+ \DD v$ of degree $ \leq k-1$. Estimate first $\deg
\widetilde \oo $ as a series. After checking divisibility we get:
$\deg \widetilde\oo = \oo \ee_0^l \leq k-1+k(T_0 -1)=k T_0 -1= \dd
q^{7 \dd+1}+ 2 \dd q^{4 \dd}$. Now using this estimation we can get
the estimation for $u,v$, it gives us  $r_{A,B} \leq \dd q^{7
\dd+1}+ 2 \dd q^{4 \dd}+ \dd$. Since $q \leq 3$ we can estimate the
latter as follows:
$r_{A,B} \leq (q+1) \dd q^{7 \dd}$. \square.

\vspace{5mm}

Let us note that using the  same idea we can get an estimation for
the degree of the generator of an infinite part of centralizer of a
given matrix $A$. We omit here details, they are similar but easier
then those were discussed above. The resulting estimation presented
in the Theorem \ref{teob}.

\section{Conjugacy separability of $\Gl(2,\BF[x])$}
\label{finsec}

In this section we would like to note that $\Gl(2,\BF[x])$ is
conjugacy separable group, but it does not immediately lead to any
algorithm which decides conjugacy in $\Gl(2,\BF[x])$, because it is
not finitely generated (see for example  \cite{Serre}) and the
finite images together with the homomorphisms onto them can not be
constructed. We mean here Maltsev's algorithm \cite{Mal} for the
decision of the conjugacy problem for finitely presented conjugacy
separable groups. Note moreover that Maltsev's algorithm does not
allow to give any estimations. To be precise let us show here the
conjugacy separability of the group $\Gl(2,\BF[x])$.

\begin{proposition}
\label{Theorem 0. }
 $\Gl(2,\BF[x])$
is conjugacy separable group.
\end{proposition}
{\bf Proof. } It is known due to Serre \cite{Serre} and Nagao
\cite{Nag} that  $\Gl(2,\BF[x]) = T(\BF[x]) \times_{T(\BF)}
\Gl(2,\BF) $ is an amalgamated free product of subgroup of upper
triangular matrices $T(\BF[x])$ and  $\Gl(2,\BF)$ through the upper
triangular matrices  over $\BF$.

There exists
the result of  J.L.Dyer  \cite{Dyer}  saying
that conjugacy separable groups
amalgamating along the finite subgroup is conjugacy separable.

We have only to verify that the subgroup $T(\BF[x])$ of upper
triangular matrices in  $\Gl(2,\BF[x])$ is conjugacy separable.

\begin{lemma}
\label{Lemma 0.}
Two elements
$ e_1=\left( \begin{array}{cc}
                \a&c\\
                0&\b
        \end {array} \right)$, \
$ e_2=\left( \begin {array} {cc}
                \a'&c'\\
                0&\b'
        \end {array} \right)$
from $T(\BF[x])$ are not conjugate if and only if

I. $(\a,\b) \neq (\a',\b') $ or

II.$(\a,\b) = (\a',\b') $ with $\a=\b$, and $c$, $c'$ are
non-proportional.
\end{lemma}

\vspace{3mm}

Let us show that for any pair $e_1 \not\sim e_2$, $e_1,e_2 \in
T(\BF[x])$ we can find normal subgroup $H_{(e_1, e_2)} \triangleleft
T(\BF[x])$, such that $\bar e_1 \not\sim  \bar e_2$, where  $\bar
e_i$ is image of $e_i$ in the finite quotient $T(\BF[x]) /H_{(e_1,
e_2)}.$

Note first that subgroups of the type
$$H_n=\left\{
             \left( \begin {array} {cc}
                1&a\\
                0&1
        \end {array} \right),
a=\a x^n+ \dots, \a \in {\BF}^* \right\}$$ are normal in
$T(\BF[x])$.

Indeed,
$\left( \begin {array} {cc}
                1&a\\
                0&1
        \end {array} \right)^e   =
\left( \begin {array} {cc}
                1&\frac{\a}{\b}a\\
                0&1
        \end {array} \right), $
where
$e=\left( \begin {array} {cc}
                \a&d\\
                0&\b
        \end {array} \right).$

To separate non-conjugate elements of the type I it is enough
to take a subgroup $H_0$.
Pick two non-conjugate elements of the type II:
$h_1=\left( \begin {array} {cc}
                \gg&c\\
                0&\gg
        \end {array} \right)$ and
$h_2=\left( \begin {array} {cc}
                \gg&c'\\
                0&\gg
        \end {array} \right).$
Note that
$
\left( \begin {array} {cc}
                \gg&c\\
                0&\gg
        \end {array} \right)^e=
\left( \begin {array} {cc}
                \gg&c'\\
                0&\gg
        \end {array} \right)
\left( \begin {array} {cc}
                1&f\\
                0&1
        \end {array} \right)
$ if and only if $\frac{\a}{\b}c+c'=\gg f$. Hence  if we take
$H_{(h_1,h_2)}= H_n$ with $n>{\rm max} \{ {\rm deg\,} c, {\rm deg\,}
c' \}$ then non-conjugates $h_1$ and $h_2$ in the quotient
$T(\BF[x])/H_{(h_1,h_2)}$ remain  non-conjugate. \square

\section{Acknowledgments}

This work has been done during the stay of the second named author in the
University of D\"usseldorf supported by DFG and in the
Max-Planck-Intit\"ut f\"ur Mathematik in Bonn, who is grateful to
these institutions for their hospitality and encouraging research
atmosphere.

 \end{document}